\documentclass{notices}
\usepackage{amsfonts,amssymb,amsmath,amscd,mathtools}

\makeatletter \@addtoreset{equation}{section}

\newtheorem{theorem}{Theorem}[section]

\newtheorem{proposition}[theorem]{Proposition}

\usepackage{color}

\newcommand{\diff}{\,\mathrm{d}}

\def\balpha{\boldsymbol{\alpha}}

\newcommand\cW{\mathcal W}

\def \E{\mathbb{E}}
\def \F{\mathbb{F}}

\def \N{\mathbb{N}}

\def \P{\mathbb{P}}

\def \R{\mathbb{R}}

\def\Ac{\mathcal{A}}

\def\Cc{\mathcal{C}}

\def\Fc{\mathcal{F}}

\def\Pc{\mathcal{P}}

\newcommand{\x}{\mathbf{x}}

\title{
On the population size in stochastic differential games
}

\author{
  Dylan Possama\"i
  \affil{
    The first author is a full professor of mathematics at ETH Z\"urich, Switzerland. 
    His email address is dylan.possamai@math.ethz.ch.
    }
  \and
  Ludovic Tangpi
  \affil{
      The second author is an assistant professor of operations research and financial engineering at Princeton University, USA. 
    His email address is ludovic.tangpi@princeton.edu.
   }
}

\begin{document}

\maketitle

Commuters looking for the shortest path to their destinations, the security of networked computers, hedge funds trading on the same stocks, governments and populations acting to mitigate an epidemic, or employers and employees agreeing on a contact, are all examples of (dynamic) stochastic differential games.
In essence, \emph{game theory} deals with the analysis of strategic interactions among multiple decision-makers.
The theory has had enormous impact in a wide variety of fields, but its rigorous mathematical analysis is rather recent.
It started with the pioneering work of von Neumann and Morgenstern \cite{vonNeum-Morg44} published in 1944. 
Since then, game theory has taken centre stage in applied mathematics and related areas, especially in economics with several game theorists such as John F. Nash Jr, Robert J. Aumann and Thomas C. Schelling being awarded the Nobel Memorial Prize in Economics Sciences.
Game theory has also played an important role in unsuspected areas: for instance in military applications, when the analysis of guided interceptor missiles in the 1950s motivated the study of games evolving dynamically in time.
Such games (when possibly subject to randomness) are called \emph{stochastic differential} games.
Their study started with the work of Issacs \cite{Issacs65}, who crucially recognised the importance of (stochastic) control theory in the area.
Over the past few decades since Isaacs' work, a rich theory of stochastic differential game has emerged and branched into several directions.
This paper will review recent advances in the study of solvability of stochastic differential games, with a focus on a purely probabilistic technique to approach the problem.
Unsurprisingly, the number of players involved in the game is a major factor of the analysis.
We will explain how the size of the population impacts the analyses and solvability of the problem, and discuss mean field games as well as the convergence of finite player games to mean field games.

\section{Two-player games}
Games involving only two players are arguably the most basic differential games in continuous time. 
In this section we discuss both zero-sum and non--zero-sum games, where the main difference stems from the existence of some symmetry---or maybe more precisely anti-symmetry here---between the players' objectives. 
Since our goal is to provide intuition, we will use a simple example as our Ariadne's thread throughout the paper, emphasising the differences and new features emerging as we make the modelling more complex.

 \smallskip
 We fix throughout this section a probability space $(\Omega,\Fc,\P)$ carrying a one-dimensional Brownian motion $W$, and we let $\F$ be the $\P$-completed natural filtration of $W$. 
 We also fix a time horizon $T>0$, and assume that $\Fc=\Fc_T$. The players are identified by numbers $1$ and $2$, and their (uncontrolled) state process is given, for some $\xi$, by
 {\small\[
 X_t\coloneqq \xi+ W_t,\; t\in[0,T].
 \]}
The players choose processes $\balpha\coloneqq (\alpha^1,\alpha^2)$, where both $\alpha^1$ and $\alpha^2$ are taken in some set of so-called \emph{admissible controls}, denoted $\Ac$, which we take here to consist essentially of functions of Brownian paths\footnote{More precisely we consider $\F$-predictable processes.} and taking values in $\R^2$ for simplicity. 
That is, players make decisions based on the randomness of the problem given by $W$.
Such controls are called \emph{open-loop}.
In fact, throughout these notes we consider only open-loop controls.
 We adopt here the so-called \emph{weak formulation}, and consider that a choice $\balpha$ from the players generates uncertainty on the distribution of $X$, by considering it under a new probability measure $\P^{\balpha}$ (it is implicitly assumed that definition of $\Ac$ ensures that this probability measure is well-defined) whose density with respect to $\P$ is given by
{\small\[
 \frac{\mathrm{d}\P^{\balpha}}{\mathrm{d}\P}\coloneqq \exp\bigg(\int_0^T(\alpha^1_s+\alpha^2_s)\diff W_s-\int_0^T\frac{(\alpha^1_s+\alpha^2_s)^2}2\diff s\bigg).
\]}
Using standard results of stochastic calculus\footnote{Notably Girsanov's theorem.}, there is then an $(\F,\P^{\balpha})$--Brownian motion $W^{\balpha}$ such that 
{\small\[
X_t\coloneqq \xi+\int_0^t(\alpha^1_s+\alpha^2_s)\diff s+ W^{\balpha}_t,\; t\in[0,T].
\]}
Of course we need to clarify how the players decide which controls they would like to use, and this will be linked to their \emph{criterion}. 
These criteria will also be exactly what is going to differentiate zero-sum and non--zero-sum games, which we will exemplify in the next sections.

\subsection{Zero-sum games and equilibria}

In zero-sum games, the goals pursued by the two players are antagonistic in the following sense:
Suppose that player's $1$ objective is to minimise to functional criterion
{\small\[
J(\alpha^1,\alpha^2)\coloneqq \E^{\P^{\balpha}}\bigg[\int_0^T\bigg(\frac{|\alpha_s^1|^2-|\alpha_s^2|^2}2+X_s\bigg)\diff s + |X_T|^2\bigg],
\]}
with $(\alpha^1,\alpha^2)\in\Ac^2$.
That is, given some control $\alpha^2\in\Ac$ chosen by player 2, player 1 aims at solving the minimisation problem\footnote{We will not discuss practical applications in these notes. To simplify the exposition, we 
rather present generic examples which extend to more general situations we will consider later. For practical examples, we refer the readers for instance to \cite{carmona2016lectures}.}
{\small\[
V^1(\alpha^2)\coloneqq \inf_{\alpha\in\Ac}J(\alpha,\alpha^2).
\]}
Then, given some control $\alpha^1\in\Ac$ chosen by player $1$, player $2$ aims at solving the minimisation problem
{\small\[
V^2(\alpha^1)\coloneqq \inf_{\alpha\in\Ac}\big\{-J(\alpha^1,\alpha^2)\big\}=-\sup_{\alpha\in\Ac}J(\alpha^1,\alpha).
\]}
In other words, while player $1$ \emph{minimises} $J$ player $2$  \emph{maximises} it.
Thus the term `zero-sum', as the sum of the criterion of each players is always $0$. 
The process $X$ is the (common) path or state of the players $1$ and $2$, $\alpha^1$ and $\alpha^2$ are their respective actions and $J(\alpha^1,\alpha^2)$ the cost (resp. reward) of player $1$ (resp. player $2$).
In this game, we are interested in so-called \emph{Nash equilibrium} corresponding to a pair $(\hat\alpha^1,\hat\alpha^2)\in\Ac^2$ such that
{\small\[
V^1(\hat\alpha^2)=J(\hat\alpha^1,\hat\alpha^2)=-V^2(\hat\alpha^1).
\]}
Thus, whenever one player plays the control corresponding to the Nash equilibrium, the other player will never be better off by not also playing according to the equilibrium. 
Observe that the Nash equilibrium is not necessarily optimal in the sense that it does not yield the highest reward (or lowest cost) to any one player, it is simply a set of strategies that is simultaneously optimal for both players.

\smallskip
The above definition of Nash equilibria anticipates our soon to come discussion of non--zero-sum games in Sect. \ref{sec:N.player}.
But this not the most standard way to attack zero-sum games. Indeed, using the (anti-)symmetry between the players' objectives, one can directly realise that if one defines instead a \emph{saddle-point} as being a process $\hat\alpha\coloneqq (\hat\alpha^1,\hat\alpha^2)\in\Ac^2$ such that
{\small\[
J(\hat\alpha^1,\hat\alpha^2)=\inf_{\alpha^ {1}\in\Ac}\sup_{\alpha^ {2}\in\Ac}J(\alpha^1,\alpha^2)=\sup_{\alpha^ {2}\in\Ac}\inf_{\alpha^ {1}\in\Ac}J(\alpha^1,\alpha^2),
\]}
then it is also a Nash equilibrium. It is important to notice at this stage that while the notion of Nash equilibrium extends to non--zero-sum games, that of saddle-points makes sense only in a zero-sum setting. Moreover, their interpretation is a bit different: a saddle-point corresponds to each player trying to optimise their criterion assuming that the other player has chosen the \emph{worst possible} control from their point of view. In general, even in the symmetric case we are describing here, not all Nash equilibria are saddle-points, but finding a saddle-point is a standard way to identify a Nash equilibrium.

\subsection{Intuition and solution}

Let us now use the simple setting we are considering to explain how one can, in general, find a saddle-point or a Nash equilibrium. 
Exactly as in standard control problems involving only one player, there are two main tools available: analytic one using partial differential equations (PDEs for short)---the celebrated Hamilton--Jacobi--Bellman (HJB for short) or Hamilton--Jacobi--Bellman--Isaacs equations---and a probabilistic one using so-called backward stochastic differential equations (BSDEs for short). This work accents the probabilistic approach (more specifically in the weak formulation), and we refer the interested reader to \cite{fleming2006controlled} for more details on the alternative one. 

\smallskip
Despite these distinctions, the point of both methods is exactly the same: understanding the structure of the so-called \emph{best-reaction function} of a player, namely what we defined as $V^1(\alpha^2)$ and $V^2(\alpha^1)$ above. Both these quantities are actually value functions of a control problem faced by each player. 
As such, it has become part of the folklore in the corresponding literature that under mild assumptions, the following result will be true. 
\begin{proposition}[Best-reaction functions]\label{prop:prop}
Under suitable assumptions, for $(\alpha^1,\alpha^2)\in\Ac^2$, we have $V^1(\alpha^2)=Y^1_0(\alpha^2)$ and $V^2(\alpha^1)=Y^2_0(\alpha^1)$ where $\big(Y^1(\alpha^2),Z^1(\alpha^2)\big)$ and $\big(Y^2(\alpha^1),Z^2(\alpha^1)\big)$ are pairs of processes satisfying appropriate measurability and integrability conditions and solving the following one dimensional \emph{BSDEs}, with $t\in[0,T]$

{\footnotesize\begin{align}
\notag
Y_t^1(\alpha^2)&=|X_T|^2+\int_t^T\bigg(\alpha^2_sZ^1_s(\alpha^2_s)-\frac{|Z^1_s(\alpha^2_s)|^2+|\alpha^2_s|^2}2\\\label{eq:bsde1}
&\quad +X_s\bigg)\mathrm{d}s-\int_t^TZ_s^1(\alpha^2_s)\mathrm{d}W_s,\\\notag
Y_t^2(\alpha^1)&=-|X_T|^2+\int_t^T\bigg(\alpha^1_sZ^2_s(\alpha^1_s)-\frac{|Z^2_s(\alpha^1_s)|^2+|\alpha^1_s|^2}2\\
\label{eq:bsde2}
&\quad -X_s\bigg)\mathrm{d}s-\int_t^TZ_s^2(\alpha^1_s)\mathrm{d}W_s.
\end{align}}
Moreover, the corresponding optimal controls are $-Z^1(\alpha^2)$ for player $1$, and $-Z^2(\alpha^1)$ for player $2$.
\end{proposition}
This result relies on appropriately using the so-called \emph{dynamic programming principle} and the associated \emph{martingale optimality principle}, as well as results from BSDE theory.
We refer for instance to \cite{zhang2017backward} for details. 
The intuition behind these equations is the following:
\begin{itemize}
\item[$(i)$] $Y_t^1(\alpha^2)$ (resp. $Y_t^2(\alpha^1)$) represents the value at time $t\in[0,T]$ of the (natural) dynamic version of the value function of player $1$ (resp. player $2$) whenever player $2$ (resp. player $1$) has chosen the control $\alpha^2$ (resp. $\alpha^1$). 
In essence, this corresponds to looking at the game over the time period $[t,T]$ only;

\smallskip
\item[$(ii)$] the functions appearing in the Lebesgue integrals in \eqref{eq:bsde1} and \eqref{eq:bsde2} are linked to the Hamiltonians of each player, which in this example are given by maps $H^1$ and $H^2$ defined on $\R^3$ by
{\small\begin{align*}
H^1(x,z,v)&\coloneqq \inf_{u\in\R}h^1(x,z,u,v)=vz-\frac{z^2+v^2}2+x,\\
H^2(x,z,u)&\coloneqq \inf_{v\in\R}h^2(x,z,u,v)=uz-\frac{z^2+u^2}2-x,
\end{align*}}
where for $(x,z,u,v)\in\R^4$,
{\small\begin{align*}
  h^1(x,z,u,v)\coloneqq (u+v)z+\frac12(u^2-v^2)+x,\\
   h^2(x,z,u,v)\coloneqq (u+a)z+\frac12(v^2-u^2)-x;
\end{align*}}

\item[$(iii)$] the processes $Z^1(\alpha^2)$ and $Z^2(\alpha^1)$ should be understood at an informal level as `derivatives'\footnote{This statement can be made rigorous using for instance Malliavin calculus.}
of $Y^1(\alpha^2)$ and $Y^2(\alpha^1)$, respectively. In practice, they directly allow to compute the optimal control for player $1$ (resp. player $2$) when player $2$ (resp. player $1$) plays $\alpha^2$ (resp. $\alpha^1)$ in the sense that it corresponds to any maximiser in the definition of $H^1(X_\cdot,Z^1_\cdot(\alpha^2),\alpha^2_\cdot)$ (resp. $H^2(X_\cdot,Z^2_\cdot(\alpha^1),\alpha^1_\cdot)$).
\end{itemize}
Once we have Proposition \ref{prop:prop} in hand, it becomes relatively straightforward to realise that to obtain a Nash equilibrium, one should be solving \emph{simultaneously} \eqref{eq:bsde1} and \eqref{eq:bsde2}, so that the behaviours of both players are concomitantly optimal. In other words, the following result holds true.
\begin{theorem}\label{th:1}
Under suitable assumptions, a pair $(\hat\alpha^1,\hat\alpha^2)\in\Ac^2$ will be a Nash equilibrium if and only if $(\hat\alpha^1,\hat\alpha^2)=(-Z^1,-Z^2)$ where the quadruplet $(Y^1,Y^2,Z^1,Z^2)$ of processes satisfies appropriate measurability and integrability conditions and solves the $2$-dimensional \emph{BSDE system}, with $t\in[0,T]$
{\footnotesize\begin{align*}
Y_t^1&=|X_T|^2+\int_t^T\bigg(Z^2_sZ^1_s -\frac{|Z^1_s|^2+|Z^2_s|^2}2+X_s\bigg)\mathrm{d}s\\
&\quad -\int_t^TZ_s^1\mathrm{d}W_s,\\
Y_t^2&=-|X_T|^2+\int_t^T\bigg(Z^1_sZ^2_s-\frac{|Z^2_s|^2+|Z^1_s|^2}2-X_s\bigg)\mathrm{d}s\\
&\quad -\int_t^TZ_s^2\mathrm{d}W_s.
\end{align*}}
\end{theorem}

The previous theorem deserves some comments, especially on how one can find Nash equilibria from the Hamiltonians of the player: the point here is that this more or less boils down to finding `fixed-points' for the vector-valued function $(H^1,H^2)$. For any $(x,z^1,z^2)\in\R^3$, what we mean by a fixed-point here is a pair $\big(u(x,z^1,z^2),v(x,z^1,z^2)\big)\in\R^2$ such that

{\footnotesize\[
\begin{cases}
  H^1\big(x,z^1,v(x,z^1,z^2)\big)=h^1\big(x,z^1,u(x,z^1,z^2),v(x,z^1,z^2)\big),\\
   H^2\big(x,z^2,u(x,z^1,z^2)\big)=h^2\big(x,z^2,u(x,z^1,z^2),v(x,z^1,z^2)\big).
\end{cases}
\]}
In our simple example, such a fixed-point is trivial to find and is uniquely given by
{\small\[
u(x,z^1,z^2)=-z^1,\; v(x,z^1,z^2)=-z^2,
\]}
which is how we identified the Nash equilibrium in \eqref{th:1}. 

\smallskip
Moreover, we would like to insist on the fact that finding a Nash equilibrium in a generic two-player game amounts to solving a 2-dimensional BSDE system. Intuitively, the dimension of the aforementioned system should increase accordingly with the number of players, and this is exactly what we will make clear in Sect. \ref{sec:N.player} below.
However, for zero-sum games (where we look for saddle-points) something interesting happens: It is enough to solve only one equation. In order to understand why, we need to introduce the so-called upper and lower values of the game, respectively denoted by $V^+$ and $V^-$, with

{\small\[
\begin{cases}
V^+\coloneqq \inf_{\alpha^ {1}\in\Ac}\sup_{\alpha^ {2}\in\Ac}J(\alpha^1,\alpha^2),\\
 V^-\coloneqq \sup_{\alpha^ {2}\in\Ac}\inf_{\alpha^ {1}\in\Ac}J(\alpha^1,\alpha^2),
 \end{cases}
\]}
as well as the upper and lower Hamiltonians, defined for $(x,z)\in\R^2$
{\small\[
\begin{cases}
  H^-(x,z)\coloneqq \sup_{v\in\R}H^1(x,z,v)=x,\\
  H^+(x,z)\coloneqq \inf_{u\in\R}\sup_{v\in\R}h^1(x,z,u,v)=x.
\end{cases}
\]}
Even if there is one major simplification in our example since $H^+$ and $H^-$ only depend on $x$, the main point to notice is rather that in general, it holds $H^+=H^-$.
This is a minimax property which constitutes what is usually refereed to as \emph{Isaacs's condition}, and is a typical necessary condition for the existence of a saddle-point. Now, in order to characterise these two values, it is useful to rely on the best-reaction functions from Proposition \ref{prop:prop}. More precisely, one can show that computing $V^+$ amounts to formally taking an infimum over $\alpha^1$ in (the opposite of) \eqref{eq:bsde2}, while computing $V^-$ amounts to formally taking a supremum over $\alpha^2$ in \eqref{eq:bsde1}. 
This means that we are naturally led to considering the pairs processes $(Y^+,Z^+)$ and $(Y^-,Z^-)$ satisfying, for $t\in[0,T]$

{\footnotesize\begin{align*}
Y_t^-&=|X_T|^2+\int_t^T\bigg(\frac{|Z^-_s|^2-|Z^-_s|^2}2+X_s\bigg)\mathrm{d}s-\int_t^TZ_s^-\mathrm{d}W_s\\
&=|X_T|^2+\int_t^TH^-(X_s,Z^-_s)\mathrm{d}s-\int_t^TZ_s^-\mathrm{d}W_s,
\end{align*}}
and
{\footnotesize\begin{align*}
Y_t^+&=|X_T|^2+\int_t^T\bigg(\frac{|Z^+_s|^2-|Z^+_s|^2}2+X_s\bigg)\mathrm{d}s-\int_t^TZ_s^+\mathrm{d}W_s\\
&=|X_T|^2+\int_t^TH_s^+(X_s,Z_s^+)\mathrm{d}s-\int_t^TZ_s^+\mathrm{d}W_s.
\end{align*}}
From there, it is immediate (by uniqueness) that $Y^+=Y^-$, $Z^+=Z^-$, and this leads to the following result.
\begin{theorem}\label{th:2}
Under suitable assumptions, we have $V^+=Y_0=V^-$ where $(Y,Z)$ solves the \emph{BSDE}
{\small\begin{align*}
Y_t&=|X_T|^2+\int_t^TX_s\mathrm{d}s-\int_t^TZ_s\mathrm{d}W_s,\; t\in[0,T].
\end{align*}}
Moreover, this equation has a unique $($and explicit$)$ solution: 
{\small\[
Y_t=T-t+X_t^2+(T-t)X_t,\; Z_t=2X_t+T-t,\; t\in[0,T],
\]}
corresponding to the saddle-point $(-Z,Z)$ for the zero-sum game.
\end{theorem}

\section{Games with an arbitrary number of players}
\label{sec:N.player}
Let us now turn our attention to the analysis of games with $N$ players, for some integer $N\ge2$.
In the two-player games discussed thus far, the main difference that we stressed was the one between zero-sum and non--zero-sum games.
With three players or more, things get much harder as the possibility of forming coalitions becomes significant, and this interesting feature has not been studied so far in the stochastic differential game literature due to the apparent difficulty of the question.
While one can consider zero-sum and cooperative versions of the $N$-player games, we will focus here on the (fully) non-cooperative case where one is interested in Nash equilibria and no coalition is formed.
In this case, for tractability of the problem, the symmetry assumption becomes essential!
We will come back to this in the final section of the article.
By symmetry here, we mean that the game will be set up in a way that players are, roughly speaking, exchangeable.
In other words, the game is exactly the same from each player's vantage point.
Of course, the reader should remark that \emph{symmetric players} does not mean \emph{independent players}, as players will still impact each other's actions and trajectories.

\smallskip

We will again discuss a probabilistic approach to the  solvability of $N$-player stochastic differential games through a simple tractable example.
Assume that the probability space $(\Omega, \Fc,\P)$ is now rich  enough to carry $N$ independent Brownian motions $(W^1,\dots, W^N)$.
Each player is identified by  an index  $i  \in  \{1,\dots,N\}$ and seeks to solve the stochastic control problem

{\footnotesize\begin{align*}
\notag
 &V^i(\alpha^{-i})\coloneqq \\ 
 &\inf_{{\alpha} {\in}{ \Ac}}\E^{\P^{{\alpha}^{{-}{i}}{\otimes}_{i} {\alpha}{,}{N}}}\bigg[\int_0^T \bigg(\frac1N\sum_{j=1}^NX^j_s+\frac12|\alpha_s|^2 \bigg)\diff s + |X^i_T|^2 \bigg],
\end{align*}}
with $X^i_t = \xi +  W^i_t$,
where 
{\small\begin{align*}
\alpha^{-i}\coloneqq  (\alpha^1,\dots,\alpha^{i-1}, \alpha^{i+1},\dots,\alpha^N),\\ 
\alpha^{-i}\otimes_i\alpha\coloneqq  (\alpha^1,\dots,\alpha^{i-1},\alpha, \alpha^{i+1},\dots,\alpha^N),
\end{align*}}
 and for a vector $\balpha$, the probability measure 
 $\P^{\balpha,N}$ is given, for any vector $\balpha\coloneqq  (\alpha^1,\dots,\alpha^N)\in\Ac^N$, by 
 
{\footnotesize\begin{equation*}
  \frac{\diff \P^{\balpha,N}}{\diff \P}\coloneqq  \exp\bigg(\sum_{j=1}^N\int_0^T\alpha^j_s\diff W_s^j-\frac12\sum_{j=1}^N\int_0^T|\alpha^j_s|^2\diff s \bigg).
\end{equation*}}

Before going any further, let us describe this problem.
Player $i$ has a state process (or trajectory) $X^i$ and control process $\alpha^i$ which is chosen among \emph{admissible} controls  $\Ac$, which are $\F^N$-predictable, where $\F^N$ is the ($\P$-completed) filtration generated by the Brownian motions $(W^1,\dots,W^N)$, and satisfy appropriate integrability conditions. This means that information on their controls is entirely encapsulated in the randomness sources $(W^1,\dots,W^N)$.
Player $i$'s goal is thus to minimise the objective function
{\small\begin{align*}
  &J(\alpha^1,\dots,\alpha^N)\coloneqq  \\
  &\E^{\P^{{\alpha}^{{-}{i}}{\otimes}_{i} {\alpha}{,}{N}}}\bigg[\int_0^T\bigg(\frac1N\sum_{j=1}^NX^j_s+\frac12|\alpha_s|^2 \bigg)\diff s + |X^i_T|^2 \bigg],
\end{align*}}
over $\alpha \in \Ac$, for $\alpha^{-i}\in\Ac^{N-1}$ fixed.
Here, the objective function is essentially the total energy of the player plus the average position of every player in the game.
The terminal cost $ |X_T^i|^2$ `forces' player $i$'s position to be close to zero at the horizon $T$.
Observe that since $(W^1,\dots,W^N)$ are independent, the state processes satisfy
{\small\begin{equation}
\label{eq:weak.solution}
  \mathrm{d}X^i_t = \alpha^i_t\diff t + \diff W^{\alpha,i}_t ,\; \P^{\balpha,N}\text{--a.s.}, \; X^i_0 = \xi,
\end{equation}}
where $W^{\alpha,i}_\cdot\coloneqq  W^i_\cdot -\int_0^\cdot \alpha^i_s\diff s$ is a $\P^{\balpha,N}$--Brownian motion.
Thus, the problem we describe aims at choosing the best  (in terms of minimising $J$) probability  space on which the state process satisfies \eqref{eq:weak.solution}.
Just as in the two player game, we will be interested in Nash equilibria, which, exactly as in the previous section, are defined as admissible strategies $(\hat\alpha^1, \dots, \hat\alpha^N)\in\Ac^N$ such  that  for every  $i \in \{1,\dots,N\}$, we have
{\small\begin{align*}
  &J(\hat\alpha^1, \dots, \hat\alpha^N) =\\
  &  \inf_{\alpha \in \Ac}\E^{\P^{\hat\alpha^{ {-} {i}}\otimes_ {i} \alpha,N}}\bigg[\int_0^T\frac12|\alpha_s|^2 + \frac1N\sum_{j=1}^NX^j_s\diff s + |X^i_T|^2 \bigg].
\end{align*}}
Again, the point here is that a Nash equilibrium corresponds to a situation where, whenever all players but one follow the equilibrium strategy, the one who does not is worse off.

\smallskip
Following the same intuition presented in the two-player case, one can rely on a system of backward SDEs
for solving and analysing Nash equilibria.
In fact, we have to following characterisation result, which directly extends Theorem \ref{th:1}.
\begin{theorem}
\label{prob.Ngame.charact}
  The $N$-player game admits a Nash equilibrium $\hat\alpha^N = (\hat\alpha^{1,N}, \dots,\hat\alpha^{N,N})$ if and only if
{\small  \begin{equation*}
    \hat\alpha^{i,N}_t = - Z^{i,i,N}_t,\; \text{\rm and}\; V^i(\hat\alpha^{-i,N}) = Y^{i,N}_0 ,\; \P\otimes\diff t\text{\rm --a.e.,}
  \end{equation*}}
   for $i\in \{1,\dots,N\}$, where $(Y^{i,N},Z^{i,j,N})_{(i,j)\in \{1,\dots,N\}^2}$ satisfy appropriate integrability conditions and solve the system
{\small  \begin{align}
  \notag
    Y^{i,N}_t &= |X^i_T|^2 + \int_t^T\bigg(\frac1N\sum_{j=1}^NX^j_s-\frac{|Z^{i,i,N}_s|^2}2 \bigg)\diff s\\
    \label{eq:Nplayer.system}
    & - \sum_{j=1}^N\int_t^TZ^{i,j,N}_s\diff W_s^{\hat\alpha^ {N},j},\; \P^{\hat\alpha^ {N},N}\text{\rm--a.s.}
  \end{align}}
  \end{theorem}
The characterisation given by Theorem \ref{prob.Ngame.charact} is important mostly because it allows to use the full force of stochastic calculus to analyse the game.
As far as solvability is concerned, it does not (necessarily) make the problem any easier because well-posedness of the system \eqref{eq:Nplayer.system} is no trivial matter in general.
The main issues here are the quadratic nonlinearity of the  drift and the fact that the system is multidimensional.
We will discuss below a very interesting problem for which Theorem \ref{prob.Ngame.charact} is particularly relevant.
Existence of $N$--Nash equilibria (and thus of the system \eqref{eq:Nplayer.system}) can also be guaranteed using compactness techniques or  analytic methods based on partial differential equations.
We refer the reader to \cite{carmona2018probabilisticI}*{Vol. I, Chapter 2} and references therein for in-depth discussion on solvability of $N$-player games.
Brushing well-posedness issues aside, the numerical simulation of equilibria (which  is  the end goal in most practical applications) is also a particularly concerning problem, especially so when the number of players $N$ is large.
Numerical simulations often reduce to approximating solutions of multidimensional systems such as \eqref{eq:Nplayer.system} or systems of $N$ HJB equations.
As is well known, most numerical simulation algorithms are subject to the so-called curse of dimensionality.
This is the fact that the performance of the algorithm plummets as the dimension increases!

\smallskip

When players are homogeneous (or symmetric) as discussed above and the interaction among the players is sufficiently weak in the sense that the influence of any given player on another one is of order $1/N$, and if the \emph{initial} configuration of the particle system is sufficiently chaotic (i.i.d.), then as $N$ increases, interaction plays less and less of a role.
One thus speaks of \emph{mean field interaction}.
Systems in mean field interactions have a long history in statistical physics (\emph{e.g.} thermodynamic models) and mathematical biology (\emph{e.g.} chemotaxis models).
In these areas, a meta theorem (which is by no means obvious) is that a symmetric particle system in weak interaction and with i.i.d. initial positions converges to an independent and identically distributed particle system whose evolution depends on the particle's probabilistic distribution.
This phenomenon is known as propagation of chaos as, in the limit, the initial chaotic configuration `propagates' to the entire path of the particles.
This is usually formulated by the fact that, if one focuses on the $k$ joint distribution at a given time $t$ of the first $k$ particles, then for $N$ sufficiently large it should approximately be given by the $k$-fold product distribution of a given particle.
Propagation of chaos was first studied by Kac \cite{kac1956foundations} and has had countless applications in pure and applied mathematics, physics, biology and more.
Limiting particle systems are usually said to be of McKean--Vlasov type.
That is, their dynamics depend on their own law.
Inspired by these ideas, Lasry and Lions \cites{lasry2006jeux} and Huang,  Malham\'e, and Caines \cite{huang2006large} proposed in 2006 a general method allowing to derive approximate Nash equilibria for large population games called \emph{mean field games}.

\section{Mean field games}

The central idea here is to consider a game in which (independent and identical) players interact through the probabilistic distribution of the entire population.
Lasry and Lions showed that solutions of these mean field games can be used to construct $N$-player strategies that are `nearly' Nash equilibria for $N$ large enough.
This reformulation  of the problem provides a decisive advantage for the simulation of equilibria, but the study of mean field games solutions (henceforth MFE for mean field equilibrium) is  an intrinsically very interesting mathematical problem.
This is particularly due to the fact that, because of the interaction through probabilistic distributions, the study of MFE heavily rests on analysis on the space of probability measures.
In fact, studying mean field games often results in the analysis 
 of PDEs as \eqref{eq:master} below written on the space $\Pc_2(\R)$  of Borel measures on $\R$ with finite second moment.
We will attempt to present a similar perspective using probabilistic arguments.
Once again, we focus on a simple example: the mean field game analogue of the game discussed in the previous section.

\smallskip
Let $\Cc([0,T],\Pc_2(\R))$ be the space of continuous functions on $[0,T]$ mapping into the set $\Pc_2(\R)$.
Let $\mu$ be an element of $\Cc([0,T],\Pc_2(\R))$.
Intuitively, $\mu_t$ should be thought of as the distribution of the (position) of the population at time $t$.
Assuming it to be known, a given (representative) player $i$ in the population will consider the stochastic control problem
{\small\begin{equation*}
  V^\mu \coloneqq  \inf_{\alpha \in \mathfrak{A}}\underbrace{\E^{\P^\alpha}\bigg[\int_0^T\bigg(\frac{|\alpha_s|^2}2 + \int_{\R}x\mu_s(\diff x)\bigg)\diff s + |X_T|^2\bigg]}_{=:J^ {\mu}(\alpha)}, 
\end{equation*}}
with $\mathrm{d}X_t = -kX_t\mathrm{d}t +  \mathrm{d}W_t^i$ and
{\small\[ \frac{\diff \P^\alpha}{\diff \P} \coloneqq  \exp\bigg(\int_0^T\alpha_s\diff W_s^i - \frac12\int_0^T|\alpha_s|^2\diff s \bigg).
\]}
Because the choice of the representative player $i$ is irrelevant (as the game is the same for all players), henceforth we will simply write $W$ instead of $W^i$.
A mean field equilibrium is a pair $(\hat\alpha, \hat\mu)\in  \mathfrak{A}\times\Cc([0,T],\Pc_2(\R))$ which satisfies
{\small\begin{equation}
\label{eq:consistency}
  V^{\hat\mu} = J^{\hat\mu}(\hat\alpha) ,\; \text{and}\; \P^{\hat\alpha}[X_t \in \cdot] = \hat\mu_t(\cdot),\; \P\otimes \diff t\text{--a.e.}
\end{equation}}
In this game, $W$ is a standard $\P$--Brownian motion. 
The set of admissible strategies $\mathfrak{A}$ is the set of processes predictable with respect to the ($\P$-completion of the) filtration generated by $W$, with appropriate integrability requirements.
The measure $\mu_t$ is a possible distribution of the path of all players  at time  $t$.
Given such a density, since all players are identical and rational, they will come up with the same optimal strategy $\hat\alpha^{\mu}$.
Therefore, at equilibrium, we expect \eqref{eq:consistency} to be satisfied.
This is precisely the same as Nash equilibrium in the sense that if all players use $\hat\alpha$, then the law of the trajectories is $\P^{\hat\alpha}[X_t\in \cdot]$, so that for a (single) player, not using $\hat\alpha$ will be suboptimal.
The condition $\P^{\hat\alpha}[X_t\in \cdot] = \hat\mu_t$ is a fixed-point condition sometimes called \emph{consistency condition}.  

\smallskip

Analogous (and actually simpler) computations leading to Theorem \ref{prob.Ngame.charact} allow to derive the following probabilistic characterisation of mean field equilibria.
\begin{theorem}
  \label{pro:charact.MFG}
  The mean field game admits a mean field equilibrium $(\hat\alpha,\hat\mu)$ if and only if it holds
  {\small\begin{equation}
  \label{eq:charaterization.MFE}
    \hat\alpha_t = - Z_t ,\; \text{\rm and}\;  V^{\hat\mu}  = Y_0,\; \P\otimes \diff t \text{\rm --a.e.},
  \end{equation}}
  where $\hat\mu_t(\cdot) = \P^{\hat\alpha}[X_t\in \cdot]$, $(Y,Z)$ satisfies appropriate integrability conditions and solves the generalised McKean--Vlasov equation
{\small  
\begin{align}
  \notag
      Y_t &= |X_T|^2 + \int_t^T\bigg( \E^{\P^{\hat\alpha}}\big[X_s]-\frac12|Z_s|^2\bigg)\diff s\\
      \label{eq:Gen.MkV.BSDE}
      &\quad - \int_t^TZ_s\diff W^{\hat\alpha}_s,\; \P^{\hat\alpha}\text{\rm--a.s.},
  \end{align}}
with $W^{\hat\alpha}_\cdot \coloneqq  W_\cdot-\int_0^\cdot\hat\alpha_s\diff s,$ and
{\small\begin{equation*}
   \frac{\diff \P^{\hat\alpha}}{\diff \P} \coloneqq  \exp\bigg(\int_0^T\hat\alpha_s\diff W_s - \frac12\int_0^T|\hat\alpha_s|^2\diff s\bigg).
  \end{equation*}}
\end{theorem}
 Equation \eqref{eq:Gen.MkV.BSDE} is a backward SDE similar to the ones encountered many times in these notes.
 It is said to be of McKean--Vlasov type because the evolution of $(Y,Z)$ depends on the law, and we call it \emph{generalised} to stress the fact that the underlying probability measure $\P^{\hat\alpha}$ as well as the underlying Brownian motion $W^{\hat\alpha}$ should be constructed together with the solution, unlike in standard BSDEs where these are given: this is a consequence of the weak formulation of the game considered here.
 For games in the strong formulation extensively discussed in \cite{carmona2018probabilisticI}, a stochastic characterisation can be derived on the fixed probability space $(\Omega, \Fc,\P)$, but the corresponding equation is then a fully-coupled system of forward--backward SDEs of McKean--Vlasov type.
 The well-posedness of \eqref{eq:charaterization.MFE} is studied in \cite{possamai2021non}.
 There exist several results on well-posedness of mean field games, see e.g. \cite{carmona2016meana},  \cite{carmona2018probabilisticI} and the references therein.
There is one extra important fact we would like to highlight regarding the above characterisation result:
\begin{itemize}
\item[$(i)$] as mentioned earlier, BSDEs characterising Nash equilibria in $N$-player games are typically of dimension $N$ themselves;

\item[$(ii)$] there is one notable exception to this rule in the case of two-player zero-sum games, since the symmetry there allows to actually consider a single equation;

\item[$(iii)$] the mean field game framework is similar, at least in spirit: thanks to the overall symmetry of the problem, we can get a characterisation using a one-dimensional equation, despite  having infinitely many players.
\end{itemize}
These observations indicate the reason why both zero-sum games and mean field games are more amenable to numerical simulations than $N$-player games, a fact that is well-recognised even for static and deterministic games.
It is also important to point out that so far, the noise in the games we have considered are only \emph{idiosyncratic} in the sense that all players deal with completely independent sources of randomness.
In most applications, a common noise is added to the problem.
In such cases, due to issues of measurability, the notion of equilibrium typically needs to be weaken or at least more carefully defined; and for solvability, the set of admissible controls sometimes needs to be enlarged to also consider relaxed (or randomised) controls  described in the next subsection.
See \cite{carmona2016meana} for details.
In addition, uniqueness of the mean field equilibrium is a delicate issue that often requires more than mere smoothness assumptions.
Structural properties and/or monotonicity properties (as the so-called Lasry and Lions monotonicity condition \cite{cardaliaguet2019master} or the displacement monotonicity condition \cite{jackson2023quantitative}) or sufficient dissipativity of the drift are usually needed.
We will refrain from discussing these---delicate---aspects of the problem here and rather direct our attention to another issue: the link between finite population and mean field games.

 \subsection{Convergence to the mean field game}
To this point, the mean field game has been motivated by heuristic arguments inspired from statistical physics 
explaining how well-behaved, uncontrolled interacting particle systems converge to independent ones in the macroscopic limit.
An interesting question is of course whether (and in which sense) stochastic differential games actually converge to mean field games.
This question is interesting in practice because it justifies the mean field game as a proxy of the $N$-player game in a large population, but it is also mathematically interesting in that it would show to which extent propagation of chaos generalises beyond uncontrolled particles systems to games.

\smallskip

This important problem was considered early on by \cites{lasry2006jeux} using PDE arguments.
The authors considered a case in which, in the $N$-player game, each player has only \emph{partial information}.
More precisely, the strategy of each player $i \in \{1,\dots,N\}$ depends only on their state process $X^i$ or their private Brownian motion $W^i$.
Under this assumption, the controlled states $X^i_t = \xi + \int_0^t-kX^i_s+\alpha^i_s\diff s + W^i_t$, becomes independent, which is a very useful simplification to prove convergence, but the partial information assumption is very restrictive in practice.
Further results have been obtained in the case of finite state space $\Omega$.
First general results are due to \cite{lacker2016general} and \cite{fischer2017connection}.
The main idea underlying the results of these authors begins with the embedding $a \longrightarrow \delta_a$ of the set $A$ in the set $\Pc(A)$ of probability measures on it.
Using this embedding, one can enlarge the set of admissible controls to consider \emph{relaxed controls}, which are measures $q$ on $[0,T]\times A$ such that $q(\diff t, \diff a) = \diff tq_t(\diff a)$ for some $q_t \in \Pc(A)$ that is sufficiently integrable.
Any (actual or strict) control $\alpha_t$ gives rise to a relaxed control $q$ given by $q(\diff t, \diff a) \coloneqq \diff t\delta_{\alpha_t}(\diff a) $.
Considering relaxed controls, \cites{lacker2016general,fischer2017connection} essentially show that the sequence of empirical laws of any Nash equilibria is relatively compact and any limiting point is a mean field equilibrium.
Further extensions can be obtained when considering closed-loop controls.

\smallskip

These general results suggest several questions; the most interesting of which, in the authors' opinion, being whether one can obtain quantitative convergence results.
That is, is it possible to obtain (non-asymptotic) convergence rates informing `how far' the $N$-player game is from the mean field game?
The simplest such bound is given by the central limit theorem which claims that, for $N$ i.i.d. square integrable random variables $(\xi^i, \dots,\xi^N)$, it holds
{\small\begin{equation*}
  \E\bigg[\bigg| \frac1N\sum_{i=1}^N\xi^i - \E[\xi^i]\bigg|^2\bigg]\le \frac{\text{Var}[\xi]}{N},\; \forall N\in\N^\star,
\end{equation*}}
where $\mathrm{Var}(\xi)$ is the variance of $\xi^1$.
Moreover, it seems natural to expect that mean field games would be explained as a form of propagation of chaos just as mean field models are for uncontrolled interacting particle systems.

\smallskip
As it turns out, Theorem \ref{prob.Ngame.charact} and Theorem \ref{pro:charact.MFG} are central in answering the above two questions.
Just as interacting particle systems model the evolution (or behaviour) of stochastically interacting components in time, the $N$-player game can also be described by the interacting particle system given by \eqref{eq:Nplayer.system}.
However, in contrast to standard particle systems, \eqref{eq:Nplayer.system} evolves \emph{backward} in time in the sense that the terminal configuration of $(Y^{i,N}, \dots, Y^{i,N})$ is known (and i.i.d.), rather than the initial one.
Similar to the mean field theory in which the mean field limit is obtained by propagation of chaos, it is now natural to expect that a form of propagation of chaos will allow to show that the (interacting) particle system defining $Y^{i,N}$ converges to $Y^i$, the i.i.d. system given by \eqref{eq:charaterization.MFE} and characterising the mean field game.
This is formalised below.
\begin{theorem}
\label{thm:main.conv}
  If for every $N\in \N^\star$ the $N$-player game admits a Nash equilibrium $(\hat\alpha^{1,N},\dots,\hat\alpha^{N,N})$ and $k$ is large enough, then the value of the $N$-player game converges to that of the mean field game in the following sense
 {\small \begin{equation*}
    |V^{i,N}(\hat\alpha^{-i,N}) - V^{\hat\mu}|^2 \le \frac{C}{N},\; \forall N\in \N^\star,
  \end{equation*}}
  for a constant $C$ depending only on $T$, and where $\hat\mu = \P^{\hat\alpha}\circ X^{-1}$.
  Moreover, for all $N\in\N^\star$
 {\small \begin{equation*}
    \E^{\P^{\hat\balpha}}\bigg[\int_0^T\cW_2^2\big(\P^{\hat\alpha^N,N}\circ (\hat\alpha^{i,N}_t)^{-1}, \P^{\hat\alpha}\circ (\hat\alpha^{i}_t)^{-1} \big)\diff t \bigg] \le \frac{C}{N}
  \end{equation*}}
  where $\cW_2$ is the Wasserstein distance.
\end{theorem}
This result is a particular case of the results derived in \cite{possamai2021non} to which we refer for detailed proofs. 
Theorem \ref{thm:main.conv} quantifies the convergence statement in the sense that it provides a rate for the convergence of the $N$-player game to the MFG.
In view of the rate dictated by the central limit theorem, this rate is sharp.
Moreover, this result provides convergence of the full sequence of values and controls, not that of a subsequence.
This result can also be established under a type of monotonicity condition known as displacement monotonicity, see \cite{jackson2023quantitative}.
In this case, we can allow much more general forms of coefficients, including controlled volatility, common noise and infinite horizon.

\subsection{Link to partial differential equations}

Most---probably all---of the results discussed in these notes have PDE counterparts that can be derived using purely analytics arguments.
Let us briefly elaborate on the link between the newly defined particle systems \eqref{eq:Nplayer.system}--\eqref{eq:Gen.MkV.BSDE} and PDE formulations.
By the so-called nonlinear Feynman--Kac formula of Peng \cite{peng1991probabilistic}, the solution $(Y^{i,N}, Z^{i,j,N})_{(i,j)\in \{1,\dots,N\}^2}$ of \eqref{eq:Nplayer.system} satisfies
{\small\begin{equation*}
  Y^{i,N}_t = v^{i,N}(t, X^{1}_t, \dots, X^{N}_t),
\end{equation*}}
where the functions $v^{i,N}:[0,T]\times \R^N\longrightarrow\R$, $i\in \{1,\dots,N\}$, is a classical solutions of the PDE system, with $\x = (x_1,\dots, x_N)$
{\small\begin{equation}
\label{eq:hjb}
  \begin{cases}
    \displaystyle\partial_tv^{i,N}(t, \x) + \frac{1}{2}\sum_{j=1}^N\partial_{x_ {j}x_ {j}}v^{i,N}(t,\x) +\frac12\big(\partial_{x_ {i}}v^{i,N}(t,\x)\big)^2\\
    \displaystyle \quad - \sum_{j=1}^N\big( \partial_{x_ {j}}v^{j,N}(t,\x)\big)^2 + \frac1N\sum_{j=1}^Nx_j = 0,\\
  \displaystyle v^{i,N}(t, \x) =  (x^i)^2,\; (t, \x) \in [0,T]\times\R^N.
  \end{cases}
\end{equation}}
This is nothing but the Hamilton--Jacobi--Bellman equation associated with the $N$-player game described earlier, see \emph{e.g.} \cite{cardaliaguet2019master}.
On the other hand, consider the so-called master equation given as 
{\small\begin{equation}
\label{eq:master}
  \begin{cases}
    \displaystyle \partial_tv(t, x, \mu) + \frac{1}{2}\partial_{xx}v(t, x,\mu) - (\partial_xv(t, x, \mu))^2\\
 \displaystyle   \quad  + \frac12 (\partial_xv(t, x, \mu))^2 - \int_{\R}\partial_xv(t, y, \mu)\partial_\mu v(t, x, \mu)(y)\\
 \displaystyle   \quad+ \frac{1}{2}\partial_x\partial_\mu v(t, x, \mu)(y)\mu(\diff y)+ \int_{\R}y\mu(\diff y)   = 0,\\
    \displaystyle v(T, x, \mu) =  x^2,\; (t,x, \mu) \in [0,T]\times \R\times \Pc_2(\R)
  \end{cases}
\end{equation}}
and characterising the mean field game (see \cite{cardaliaguet2019master}).
Here, $\partial_\mu v$ denote the $L$-derivative of the function $v$ which is understood as follows:

\smallskip
$(i)$ lift the function $v(\cdot,\cdot,\mu)$ on $L^2(\Omega,\Fc,\P^{\hat\alpha})$ by putting $\tilde v(\cdot,\cdot, \xi)\coloneqq v(\cdot,\cdot, \mu)$ whenever $\xi \in L^2(\Omega,\Fc,\P^{\hat\alpha})$ has law $\mu$; 

\smallskip
$(ii)$ denote by $\partial_\mu v(\cdot,\cdot, \mu)$ the Fr\'echet derivative of the lift $\tilde v(\cdot,\cdot, \xi)$ on $L^2(\Omega,\Fc, \P^{\hat\alpha})$.

\smallskip
\noindent Equation \eqref{eq:master} is a one-dimensional equation written on the infinite dimensional space $[0,T]\times \R\times \Pc_2(\R)$ whereas \eqref{eq:hjb} is a multi-dimensional equation written on the finite dimensional space $[0,T]\times \R^N$.
When \eqref{eq:master} admits a classical solution, it follows (see e.g. \cite{chassagneux2014probabilistic} applied on the probability space $(\Omega, \Fc, \P^{\hat\alpha})$) that
{\small\begin{equation*}
  Y_t = v(t, X_t, \P^{\hat\alpha}\circ X_t^{-1}).
\end{equation*}}
Well-posedness of \eqref{eq:master} in the classical sense is unfortunately very difficult to obtain in most cases.
The work \cite{cardaliaguet2019master} exhibits cases of games for which the master equation admits a unique solution which in addition has bounded derivatives.
In particular, as a consequence of Theorem \ref{thm:main.conv} one infers that for $i\in \{1,\dots,N\}$
{\small\begin{equation}
\label{eq:V.conv.N}
  |v^{i, N}(0, X_0,\dots, X_0) - v(0, X_0, \delta_{X_0})|^2 \le \frac CN  ,\; \forall N\in \N^\star.
\end{equation}}
\cite{cardaliaguet2019master} shows that when adequate monotonicity is satisfied, the dissipativity coefficient $k$ is not needed to guarantee \eqref{eq:V.conv.N}.

\section{Going further}

Our discussion of stochastic differential games is just a glimpse of a rich and fast developing theory.
These notes focused on how the analysis changes with the number of players, and we presented only games in the probabilistic weak formulation.
This exposition omits important issues such as general conditions guaranteeing well-posedness and properties of equilibria, or the importance of the type of information used by players.
For instance, we did not discuss the so-called closed-loop, or Markovian controls, nor the various type of monotonicity conditions on the data often needed to guarantee uniqueness of games such as the Lasry--Lions monotonicity and the displacement monotonicity. We also completely omitted the deep connections and similarities between optimal transport and mean field games: indeed, some particular mean field game problems, such as potential mean field games, can be written as a dynamical optimal transport problem. This opens the way for applying techniques from one field to the other, see for instance \cite{li2023controlling}. Moreover, we would not want our readers to think that the theory is limited to continuous-time and/or uncountable state spaces: there is plethora of contributions analysing finite state mean field games, and which contribute deeply to the theory, see among many others \cites{gomes2013continuous,Bayraktar21}.

\smallskip

Many of the `principles' derived for competitive games in this article extend, albeit with additional technicalities to other types of games which we could not discuss here.
This is for instance the case for cooperative games in which a social planner finds the best strategies for individual players to achieve a common objective.
This is the kind of problems solved by ride sharing apps.
Games with leaders and followers where some players react to other players actions or decisions (e.g. Stackelberg games) or major-minor players can also be recast in a setting similar to the one discussed here.
Such games are central in contract theory applications.
Speaking of applications, as explained in the previous section, mean field games arise from finite player games in which players are symmetric and interact with (almost) all other players in the game.
This is certainly not true for most applications.
It would be interesting to study games in which players (possibly) interact differently with each of their peers; and where the form of interaction can evolve over time.
Such games are called \emph{games on graphs}.
Although recent progress have been made, the study of games on graphs is still in its infancy, and it will undoubtedly provide very powerful modelling tools.

\smallskip
The vast majority of works on stochastic differential games, just as the present one, consider Nash equilibria.
This has allowed to develop an impressive and successful theory.
Notwithstanding this success, it is worth pointing out that the concept of Nash equilibrium itself presents several shortcomings.
For instance, Nash equilibria are rarely unique and among the many equilibria, there can be several trivial or non desirable ones, which poses the questions of \emph{selection}, \emph{approximation} and \emph{stability} of equilibria. 
In fact, unless players can negotiate to agree on an equilibrium, they can choose different equilibria, which may in turn result in a non-equilibrium set of strategies.
Even when the Nash equilibrium is unique, the rationality assumption might be too restrictive since in most practical cases agents would rather settle for a satisfactory outcome than an optimal one.
Analysing stochastic differential games beyond Nash equilibrium seems to be an interesting avenue for future research.

\bibliography{ExampleRefs}


\end{document}